\documentclass[12pt]{article}
\usepackage{amsfonts}
\usepackage{amsfonts}
\usepackage{mathrsfs}
\usepackage{amsfonts}
\usepackage{amssymb,amsfonts,amsmath, psfrag,eepic,colordvi,epsfig}
\topskip 
\numberwithin{equation}{section}
\newtheorem{theorem}{Theorem}[section]

\newtheorem{conjecture}[theorem]{Conjecture}

\newtheorem{lemma}[theorem]{Lemma}

\begin{document}
\parskip 7pt

\pagenumbering{arabic}
\def\sof{\hfill\rule{2mm}{2mm}}
\def\ls{\leq}
\def\gs{\geq}
\def\SS{\mathcal S}
\def\qq{{\bold q}}
\def\MM{\mathcal M}
\def\TT{\mathcal T}
\def\EE{\mathcal E}
\def\lsp{\mbox{lsp}}
\def\rsp{\mbox{rsp}}
\def\pf{\noindent {\it Proof.} }
\def\mp{\mbox{pyramid}}
\def\mb{\mbox{block}}
\def\mc{\mbox{cross}}
\def\qed{\hfill \rule{4pt}{7pt}}
\def\pf{\noindent {\it Proof.} }
\textheight=22cm

{\Large
\begin{center}
 Proof of a Conjecture of  Hirschhorn and Sellers
  on Overpartitions
\end{center}
}

\begin{center}
 William Y.C. Chen$^1$ and  Ernest X.W. Xia$^2$

 $^1$Center for Applied Mathematics\\
Tianjin University\\
Tianjin 300072, P.R. China

$^2$Department of Mathematics\\
 Jiangsu University\\
 Zhenjiang, Jiangsu 212013, P.R. China

 $^1$chenyc@tju.edu.cn, $^2$ernestxwxia@163.com

\end{center}

\noindent {\bf Abstract.} Let $\bar{p}(n)$ denote the number of
overpartitions of $n$. It was conjectured by   Hirschhorn and
Sellers that  $\bar{p}(40n+35)\equiv 0\ ({\rm mod\ } 40)$ for $n\geq
0$. Employing $2$-dissection formulas  of
  quotients of   theta functions due to Ramanujan, and  Hirschhorn
and Sellers,  we obtain  a generating function for $\bar{p}(40n+35)$
modulo 5. Using the $(p, k)$-parametrization of theta functions
given by Alaca, Alaca and Williams,  we give a proof  of the
congruence
  $\bar{p}(40n+35)\equiv
0\ ({\rm mod\ } 5)$. Combining this congruence and the congruence
 $\bar{p}(4n+3)\equiv
0\ ({\rm mod\ } 8)$ obtained
  by Hirschhorn and Sellers, and   Fortin,   Jacob and
    Mathieu, we give a proof of the conjecture of Hirschhorn and Sellers.

   \noindent {\bf Keywords:}  overpartition, congruence,
   theta function, $(p,k)$-parametrization.

\noindent {\bf AMS Subject
 Classification:} 11P83, 05A17.

\section{Introduction}

\allowdisplaybreaks

The objective of this paper is to give a proof of a conjecture of
Hirschhorn and Sellers on the number of overpartitons. We shall use
the technique of dissections of quotients of theta functions.

Let us begin with some notation and terminology on $q$-series and
partitions.
 We adopt the common notation
 \begin{align}
(a;q)_\infty=\prod_{n=0}^{\infty}(1-aq^n), \label{1-1}
 \end{align}
 where $|q|<1$,
 and we write
 \begin{align}
(a_1,a_2,\ldots, a_n;q)_\infty &=(a_1;q)_\infty(a_2;q)_\infty \cdots
(a_n;q)_\infty.\label{1-2}
 \end{align}
 Recall that the Ramanujan
  theta function $f(a,b)$ is defined by
    \begin{align}\label{1-3}
f(a,b)=\sum_{n=-\infty}^\infty
 a^{n(n+1)/2}b^{n(n-1)/2},
  \end{align}
where $|ab|<1$. The  Jacobi
  triple product identity can be restated as
\begin{align}\label{1-4}
f(a,b)=(-a,-b,ab;ab)_\infty.
\end{align}
Here is a  special case   of \eqref{1-3}, namely,
\begin{align}
 f(-q)=f(-q,-q^2)
 =\sum_{n=-\infty}^\infty
 (-1)^nq^{n(3n-1)/2}=(q;q)_\infty.\label{1-7}
\end{align}
For any  positive integer $n$, we use $f_n$ to denote  $f(-q^n)$,
that is,
\begin{align}
f_n=(q^n;q^n)_\infty=\prod_{k=1}^\infty (1-q^{nk}).\label{R-2}
\end{align}

The  function $f_n$ is related to the generating function of
overpartitions.  A partition of a positive integer $n$ is
 any nonincreasing sequence of positive integers whose
sum is $n$.  An overpartition of $n$ is a partition in which the
 first occurrence of a number may be overlined, see Corteel and
 Lovejoy \cite{Corteel}.
 For $n\geq 1$, let $\bar{p}(n)$
denote the number of overpartitions of
 $n$, and we set $\bar{p}(0)=1$.
 Corteel and Lovejoy \cite{Corteel} showed that
   the generating function for
$\bar{p}(n)$ is given by
\begin{align}
\sum_{n=0}^\infty \bar{p}(n)q^n=\frac{f_2}{f_1^2}.  \label{1-11}
\end{align}

 Hirschhorn and   Sellers \cite{Hirschhorn-1},
   and   Fortin,   Jacob and
  Mathieu  \cite{Fortin}
obtained the following Ramanujan-type generating function formulas
for $\bar{p}(2n+1)$, $\bar{p}(4n+3)$,
  and  $\bar{p}(8n+7)$:
\begin{align}
\sum_{n=0}^\infty \bar{p}(2n+1)q^n&= 2\frac{f_2^2f_8^2}{f_1^4f_4},
\label{M-0} \\[6pt]
\sum_{n=0}^\infty \bar{p}(4n+3)q^n&= 8\frac{f_2f_4^6}{f_1^8},
\label{M-1}
\\[6pt]
\sum_{n=0}^\infty \bar{p}(8n+7)q^n&= 64\frac{f_2^{22}}{f_1^{23}}.
\label{M-2}
\end{align}
The above identities lead to congruences modulo  2, 8 and  64 for
the overpartition function.  Mahlburg \cite{Mahlburg} proved that
$\bar{p}(n)$
     is divisible by 64 for almost all $n$ by using   relations between $\bar{p}(n)$ and the number
of representations of $n$ as a sum of squares. Using the theory of
modular forms, Treneer  \cite{Treneer}  showed that the coefficients
of a wide class of weakly holomorphic modular forms have infinitely
many congruence relations for powers of every prime $p$ other than 2
and 3. In particular,  Treneer  \cite{Treneer}  proved that
$\bar{p}(5m^3n)\equiv  0\ ({\rm mod}\  5)$ for any  $n$ that is
coprime to $m$,
 where $m$ is a prime satisfying $m \equiv -1 \ ({\rm mod}\  5)$.

The following conjecture was posed by  Hirschhorn and Sellers
\cite{Hirschhorn-1}.

\begin{conjecture} \label{conjecture-1}
For $n\geq 0$, we have
\begin{align}
\bar{p}(40n+35)\equiv 0 \quad ({\rm mod \ 40}).\label{1-12}
\end{align}
\end{conjecture}

This paper is organized as follows. In Section 2, using
$2$-dissection formulas  of  quotients of theta functions given by
Ramanujan \cite{Berndt1991}, and Hirschhorn and Sellers
\cite{Hirschhorn-3},
  we  derive
    a generating function for $\bar{p}(40n+35)$
modulo 5.   In Section 3,  we use the $(p, k)$-parametrization of
theta functions due to  Alaca, Alaca and Williams
\cite{Alaca,Alaca-1,Williams} to show that $\bar{p}(40n+35)\equiv 0\
({\rm mod\ } 5)$. This proves the conjecture of Hirschhorn and
Sellers resorting to the congruence
 $\bar{p}(4n+3)\equiv
0\ ({\rm mod\ } 8)$ independently obtained
  by Hirschhorn and Sellers \cite{Hirschhorn-1},
  and  Fortin,   Jacob and   Mathieu \cite{Fortin}.

\section{Generating function of $\bar{p}(40n+35)$ modulo  5}

In this Section, we deduce a generating function of
$\bar{p}(40n+35)$ modulo  5.
 We first recall several
  $2$-dissection formulas for quotients of theta functions
   due to   Ramanujan \cite{Berndt1991}, Hirschhorn
and Sellers \cite{Hirschhorn-3}.

The following relations are  consequences of dissection formulas of
Ramanujan collected in  Entry 25  in Berndt's book \cite[p.
40]{Berndt1991}.

\begin{lemma}\label{L-1} Let $f_n$ be given by  \eqref{R-2}.
 We have  \begin{align}
f_1^2&=\frac{f_2f_8^5}{f_4^2f_{16}^2} -2q\frac{f_2f_{16}^2}{f_8},
\label{2-1}\\[6pt]
\frac{1}{f_1^2}&=
\frac{f_8^5}{f_2^5f_{16}^2}+2q\frac{f_4^2f_{16}^2}{f_2^5f_8}, \label{2-2}\\[6pt]
f_1^4&=\frac{f_4^{10}}{f_2^2f_8^4}-4q\frac{f_2^2f_8^4}{f_4^2}
\label{2-3}
\end{align}
and
\begin{align}\label{2-4}
\frac{1}{f_1^4}&=\frac{f_4^{14}}{f_2^{14}f_8^4}
+4q\frac{f_4^2f_8^4}{f_2^{10}}.
\end{align}
\end{lemma}

Hirschhorn and   Sellers   \cite{Hirschhorn-3} established  the
following $2$-dissection formula.

 \begin{lemma}\label{L-1-1} Let $f_n$ be given by  \eqref{R-2}. We have
\begin{align}\label{2-r-1}
\frac{f_5}{f_1}=\frac{f_8f_{20}^2}{f_2^2f_{40}}
+q\frac{f_4^3f_{10}f_{40}}{f_2^3f_8f_{20}}.
\end{align}
\end{lemma}

By Lemmas \ref{L-1} and
  \ref{L-1-1}, we are led to a  generating function
  of $\bar{p}(40n+35)$ modulo 5.

\begin{theorem}\label{Th-0} We have
\begin{align}
\sum_{n=0}^\infty \bar{p}(40n+35)q^n \equiv &
2\frac{f_2^{122}}{f_1^{63}f_4^{40}} +3 \frac{f_1f_2^{26}}{f_4^8}
 +4q\frac{f_2^{98}}{f_1^{55}f_4^{24}}
 +3qf_1^9f_2^2f_4^8  +4q^2\frac{f_2^{74}}{f_1^{47}f_4^8}\nonumber\\[7pt]
&\quad
 +4q^3\frac{f_2^{50}f_4^8}{f_1^{39}}
  +4q^4\frac{f_2^{26}f_4^{24}}{f_1^{31}}  + 2q^{5}
\frac{f_2^2f_4^{40}}{f_1^{23}}\quad ({\rm mod \ 5}).\label{R-1}
\end{align}
\end{theorem}

\noindent{\it Proof.} Recall that the well-known theta function
$\varphi(q)$  is  defined by
\begin{align}
\varphi(q)&=f(q,q)=\sum_{n=-\infty}^\infty q^{n^2}.\label{1-5}
\end{align}
By the Jacobi triple product identity, we find
  \begin{align}
\varphi(q)=\frac{f_2^5}{f_1^2f_4^2}.\label{1-8-0}
\end{align}
In view  of  \eqref{1-3},   \eqref{1-4} and \eqref{1-5}, we see that
\begin{align}
\varphi(q)=\sum_{n=-\infty }^\infty q^{n^2}&=\sum_{n=-\infty}^\infty
q^{25n^2}+2q\sum_{n=-\infty}^\infty q^{25n^2+10n}
+2q^4\sum_{n=-\infty}^\infty q^{25n^2+20n}\nonumber\\[6pt]
&=\varphi(q^{25})+2qD(q^5)+2q^4E(q^5),\label{3-1}
\end{align}
where
\begin{align}
D(q)=\sum_{n=-\infty}^\infty
q^{5n^2+2n}=(-q^3,-q^7,q^{10};q^{10})_\infty \label{3-2}
\end{align}
and
\begin{align}
E(q)=\sum_{n=-\infty}^\infty q^{5n^2+4n}=
 (-q,-q^9,q^{10};q^{10})_\infty. \label{3-3}
\end{align}
It is easily checked that
\begin{align}
D(q)E(q)=\frac{f_2^2f_5f_{20}}{f_1f_4 }. \label{3-4}
\end{align}
Replacing $q$ by $-q$ in \eqref{1-8-0},  and using the fact that
\begin{align}
f(q)=(-q;-q)_\infty =\frac{f_2^3}{f_1f_4},\label{1-9}
\end{align}
we deduce that
\begin{align}
\varphi(-q)=\frac{f_1^2}{f_2} .\label{1-10-0}
\end{align}
Because of \eqref{1-10-0}, the generating function \eqref{1-11} of
$\bar{p}(n)$ can be rewritten as
\begin{align}
\sum_{n=0}^\infty \bar{p}(n)(-q)^n
=\frac{1}{\varphi(q)}.\label{R-11}
\end{align}
 It follows that
\begin{align}
\sum_{n=0}^\infty \bar{p}(n)(-q)^n
=\frac{\varphi^4(q)}{\varphi^5(q)}. \label{R-9}
\end{align}
By the binomial theorem, it is easily seen that for $k\geq 1$,
\begin{align}
(1-q^k)^{5}\equiv (1-q^{5k}) \quad   ({\rm mod} \ 5),\label{R-3}
\end{align}
which implies that
\begin{align} \label{3-5}
\varphi^5(q)\equiv \varphi(q^5)\quad ({\rm mod} \ 5).
\end{align}
Combining   \eqref{R-9} and \eqref{3-5}, we find that
\begin{align}
\sum_{n=0}^\infty \bar{p}(n)(-q)^n \equiv \frac{\varphi^4(q)}{\varphi(q^5)} \quad ({\rm mod }\ 5).\label{R-10}
\end{align}
Substituting \eqref{3-1} into \eqref{R-10}, we get
\begin{align}
\sum_{n=0}^\infty \bar{p}(n)(-q)^n
&\equiv
\frac{\left(\varphi(q^{25})+2qD(q^5)+2q^4E(q^5)
\right)^4}{\varphi(q^5)} \nonumber\\[6pt]
&\equiv \frac{1}{\varphi(q^5)}(\varphi^4(q^{25})
+3q\varphi^3(q^{25})D(q^5)+4q^2\varphi^2(q^{25})
D^2(q^5)+2q^3\varphi(q^{25})D^3(q^5)
\nonumber\\[6pt]
&\qquad +3q^4\varphi^3(q^{25})E(q^5)+q^4D^4(q^5)
+3q^5\varphi^2(q^{25})D(q^5)D(q^5)
\nonumber\\[6pt]
&\qquad +q^6\varphi(q^{25})D^2(q^5)E(q^5)
+4q^7D^3(q^5)E(q^5)+4q^8\varphi^2(q^{25})E^2(q^5)
\nonumber\\[6pt]
&\qquad +q^9\varphi(q^{25})D(q^5)E^2(q^5) +q^{10}D^2(q^5)E^2(q^5)
+2q^{12}\varphi(q^{25})E^3(q^5)\nonumber\\[6pt]
&\qquad +4q^{13}D(q^5)E^3(q^5) + q^{16}E^4(q^5)^4) \quad ({\rm mod}
\ 5).\label{3-6}
\end{align}
Extracting those  terms  associated with powers $q^{5n}$
  on both sides of
\eqref{3-6} and replacing $q^5$ by $q$, we obtain
\begin{align}
\sum_{n=0}^\infty \bar{p}(5n)(-q)^n &\equiv \frac{\varphi^4(q^{5})
+3q\varphi^2(q^{5})D(q)E(q)
  +q^{2}D^2(q)E^2(q)}{\varphi(q)} \quad ({\rm mod} \ 5).\label{3-7}
\end{align}
By \eqref{1-8-0} and \eqref{3-4}, we can rewrite  \eqref{3-7} as
follows
\begin{align}
\sum_{n=0}^\infty \bar{p}(5n)(-q)^n &\equiv
\frac{f_1^2f_4^2f_{10}^{20}}{f_2^5f_5^8f_{20}^8}
+3q\frac{f_1f_4f_{10}^{10}}{f_2^3f_5^3f_{20}^3}
+q^2\frac{f_5^2f_{20}^2}{f_2} \qquad ({\rm mod} \ 5). \label{3-8}
\end{align}
Replacing $q$ by $-q$  in  \eqref{3-8}, we get
\begin{align}
\sum_{n=0}^\infty \bar{p}(5n)q^n &\equiv
\frac{f_4^2f_{10}^{20}}{f_2^5f_{20}^8}\frac{f^2(q)}{f^8(q^5)}
-3q\frac{f_4f_{10}^{10}}{f_2^3f_{20}^3}\frac{f(q)}{f^3(q^5)}
+q^2\frac{f_{20}^2}{f_2}f^2(q^5) \qquad ({\rm mod} \ 5).
\label{R-12}
\end{align}
Substituting  \eqref{1-9} into \eqref{R-12},
 we arrive at
\begin{align}
\sum_{n=0}^\infty \bar{p}(5n)q^n &\equiv
\frac{f_2f_5^8}{f_1^2f_{10}^4}
-3q\frac{f_5^3f_{10}}{f_1}+q^2\frac{f_{10}^6}{f_2f_5^2} \quad ({\rm
mod} \ 5).\label{3-9}
\end{align}
According to   $2$-dissection formulas  \eqref{2-1}, \eqref{2-2},
\eqref{2-3}, \eqref{2-r-1} and congruence \eqref{3-9},
  we
obtain
\begin{align}
\sum_{n=0}^\infty \bar{p}(5n)q^n \equiv& \frac{f_2}{f_{10}^4}
\left(\frac{f_8^5}{f_2^5f_{16}^2}+2q\frac{f_4^2f_{16}^2}
{f_2^5f_8}\right) \left( \frac{f_{20}^{10}}{f_{10}^2f_{40}^4}
-4q^5\frac{f_{10}^2f_{40}^4}{f_{20}^2}\right)^2\nonumber\\[6pt]
&\quad  -3q f_{10}\left(\frac{f_{10}f_{40}^5}{f_{20}^2f_{80}^2}
-2q^5\frac{f_{10}f_{80}^2}{f_{40}}\right) \left(
\frac{f_8f_{20}^2}{f_2^2f_{40}}
+q\frac{f_4^3f_{10}f_{40}}{f_2^3f_8f_{20}}\right) \nonumber\\[6pt]
&\quad
+q^2\frac{f_{10}^6}{f_2}\left(\frac{f_{40}^5}{f_{10}^5f_{80}^2}
+2q^5\frac{f_{20}^2f_{80}^2} {f_{10}^5f_{40}}\right) \nonumber\\[6pt]
\equiv& \frac{f_8^5f_{20}^{20}}{f_2^4f_{10}^8f_{16}^2f_{40}^8}
+2q\frac{f_4^2f_{16}^2f_{20}^{20}}{f_2^4f_8f_{10}^8f_{40}^8}
-3q\frac{f_8f_{10}^2f_{40}^4}{f_2^2f_{80}^2}
-3q^2\frac{f_4^3f_{10}^3f_{40}^6}{f_2^3f_8f_{20}^3f_{80}^2}
\nonumber\\[6pt]
&\quad  +q^2\frac{f_{10}f_{40}^5}{f_2f_{80}^2}
-3q^5\frac{f_8^5f_{20}^8}{f_2^4f_{10}^4f_{16}^2}
+q^6\frac{f_8f_{10}^2f_{20}^2f_{80}^2}{f_2^2f_{40}^2}
-q^6\frac{f_4^2f_{16}^2f_{20}^8}{f_2^4f_8f_{10}^4}+q^7\frac{f_4^3f_{10}^3f_{80}^2}{f_2^3f_8f_{20}} \nonumber\\[6pt]
&\quad +2q^7\frac{f_{10}f_{20}^2f_{80}^2}{f_2f_{40}} + q^{10}
\frac{f_8^5f_{40}^8}{f_2^4f_{16}^2f_{20}^4} +2q^{11}
\frac{f_4^2f_{16}^2f_{40}^8}{f_2^4f_8f_{20}^4} \quad ({\rm mod} \
5).\label{3-10}
\end{align}
Extracting the terms with odd powers of $q$ on both sides of
\eqref{3-10}, we have
\begin{align}
\sum_{n=0}^\infty \bar{p}(10n+5)q^{2n+1}  \equiv&
2q\frac{f_4^2f_{16}^2f_{20}^{20}}{f_2^4f_8f_{10}^8f_{40}^8}
-3q\frac{f_8f_{10}^2f_{40}^4}{f_2^2f_{80}^2}
-3q^5\frac{f_8^5f_{20}^8}{f_2^4f_{10}^4f_{16}^2}
+q^7\frac{f_4^3f_{10}^3f_{80}^2}{f_2^3f_8f_{20}} \nonumber\\[6pt]
&\quad +2q^7\frac{f_{10}f_{20}^2f_{80}^2}{f_2f_{40}}  +2q^{11}
\frac{f_4^2f_{16}^2f_{40}^8}{f_2^4f_8f_{20}^4} \quad ({\rm mod} \
5).\label{3-10-1}
\end{align}
Dividing $q$ on both sides of \eqref{3-10-1} and   replacing $q^2$
by $q$, we get
\begin{align}
\sum_{n=0}^\infty \bar{p}(10n+5)q^n \equiv&
2\frac{f_2^2f_{8}^2f_{10}^{20}}{f_1^4f_4f_{5}^8f_{20}^8}
-3\frac{f_4f_{5}^2f_{20}^4}{f_1^2f_{40}^2}
-3q^2\frac{f_4^5f_{10}^8}{f_1^4f_{5}^4f_{8}^2}
\nonumber\\[6pt]
&\quad
 +q^3\frac{f_2^3f_{5}^3f_{40}^2}{f_1^3f_4f_{10}}  +2q^3\frac{f_{5}f_{10}^2f_{40}^2}{f_1f_{20}}  +2q^{5}
\frac{f_2^2f_{8}^2f_{20}^8}{f_1^4f_4f_{10}^4} \quad ({\rm mod} \ 5).
\label{3-11}
\end{align}
Employing $2$-dissection formulas  \eqref{2-4},  \eqref{2-r-1} and
 congruence \eqref{3-11}, we deduce that
\begin{align}
\sum_{n=0}^\infty \bar{p}(10n+5)q^n \equiv&
2\frac{f_2^2f_{8}^2f_{10}^{20}}{ f_4f_{5}^8f_{20}^8}
\left(\frac{f_4^{14}}{f_2^{14}f_8^4}
+4q\frac{f_4^2f_8^4}{f_2^{10}}\right)
\left(\frac{f_{20}^{14}}{f_{10}^{14}f_{40}^4}
+4q^5\frac{f_{20}^2f_{40}^4}{f_{10}^{10}}\right)^2\nonumber\\[6pt]
&\quad  -3\frac{f_4 f_{20}^4}{ f_{40}^2}\left(
\frac{f_8f_{20}^2}{f_2^2f_{40}}
+q\frac{f_4^3f_{10}f_{40}}{f_2^3f_8f_{20}}\right)^2 \nonumber\\[6pt]
&\quad -3q^2\frac{f_4^5f_{10}^8}{
f_{8}^2}\left(\frac{f_4^{14}}{f_2^{14}f_8^4}
+4q\frac{f_4^2f_8^4}{f_2^{10}}\right)
\left(\frac{f_{20}^{14}}{f_{10}^{14}f_{40}^4}
+4q^5\frac{f_{20}^2f_{40}^4}{f_{10}^{10}}\right)
\nonumber\\[6pt]
&\quad
 +q^3\frac{f_2^3 f_{40}^2}{ f_4f_{10}}
 \left(
\frac{f_8f_{20}^2}{f_2^2f_{40}}
+q\frac{f_4^3f_{10}f_{40}}{f_2^3f_8f_{20}}\right)^3 +2q^3\frac{
f_{10}^2f_{40}^2}{ f_{20}}\left( \frac{f_8f_{20}^2}{f_2^2f_{40}}
+q\frac{f_4^3f_{10}f_{40}}{f_2^3f_8f_{20}}\right)^2 \nonumber\\[6pt]
&\quad  +2q^{5} \frac{f_2^2f_{8}^2f_{20}^8}{
f_4f_{10}^4}\left(\frac{f_4^{14}}{f_2^{14}f_8^4}
+4q\frac{f_4^2f_8^4}{f_2^{10}}\right)\nonumber\\[6pt]
\equiv & -3\frac{f_4f_8^2f_{20}^8}{f_2^4f_{40}^4}
+2\frac{f_4^{13}f_{20}^{20}}{f_2^{12}f_8^2f_{10}^8f_{40}^8}
+3q\frac{f_4f_8^6f_{20}^{20}}{f_2^8f_{10}^8f_{40}^8}
-q\frac{f_4^4f_{10}f_{20}^5}{f_2^5 f_{40}^2} \nonumber\\[6pt]
&\quad
-3q^2\frac{f_4^{19}f_{20}^{14}}{f_2^{14}f_8^6f_{10}^6f_{40}^4}
-3q^2\frac{f_4^7f_{10}^2f_{20}^2}{f_2^6f_8^2}
-2q^3\frac{f_4^7f_8^2f_{20}^{14}}{f_2^{10}f_{10}^6f_{40}^4}
+2q^3\frac{f_8f_{10}^2f_{20}f_{40}}{f_2^2} \nonumber\\[6pt]
&\quad +q^3\frac{f_8^3f_{20}^6}{f_2^3f_4f_{10}f_{40}}
+2q^4\frac{f_4^3f_{10}^3f_{40}^3}{f_2^3f_8f_{20}^2}
+3q^4\frac{f_4^2f_8f_{20}^3f_{40}}{f_2^4}
+3q^5\frac{f_4^5f_{10}f_{40}^3}{f_2^5f_8} \nonumber\\[6pt]
&\quad +3q^5\frac{f_4^{13}f_{20}^8}{f_2^{12}f_8^2f_{10}^4}
+2q^6\frac{f_4f_8^6f_{20}^8}{f_2^8f_{10}^4}
+q^6\frac{f_4^8f_{10}^2f_{40}^5}{f_2^6f_8^3f_{20}^3}
-2q^7\frac{f_4^{19}f_{20}^2f_{40}^4}{f_2^{14}f_8^6f_{10}^2}
 \nonumber\\[6pt]
 &\quad -3q^8\frac{f_4^7f_8^2f_{20}^2f_{40}^4}{f_2^{10}f_{10}^2}
 +2q^{10}\frac{f_4^{13}f_{40}^8}{f_2^{12}f_8^2 f_{20}^4}
 +3q^{11}\frac{f_4f_8^6f_{40}^8}{f_2^8f_{20}^4} \quad ({\rm mod} \
 5).
 \label{3-12}
\end{align}
Extracting the terms with odd powers of $q$ on both sides of
\eqref{3-12}, then dividing by $q$ and replacing $q^2$ by $q$,
 we find that
\begin{align}
\sum_{n=0}^\infty \bar{p}(20n+15)q^n \equiv&
3\frac{f_2f_4^6f_{10}^{20}}{f_1^8f_5^8f_{20}^8}
-\frac{f_2^4f_5f_{10}^5}{f_1^5f_{20}^2}
-2q\frac{f_2^7f_4^2f_{10}^{14}}{f_1^{10}f_5^6f_{20}^4}
+2q\frac{f_4f_5^2f_{10}f_{20}}{f_1^2}
\nonumber\\[6pt]
&\quad  +q\frac{f_4^3f_{10}^6}{f_1^3f_2f_5f_{20}}
+3q^2\frac{f_2^5f_5f_{20}^3}{f_1^5f_4}
+3q^2\frac{f_2^{13}f_{10}^8}{f_1^{12}f_4^2f_5^4} \nonumber\\[6pt]
&\quad -2q^3\frac{f_2^{19}f_{10}^2f_{20}^4}{f_1^{14}f_4^6f_5^2}
+3q^5\frac{f_2f_4^6f_{20}^8}{f_1^8f_{10}^4}\quad   ({\rm mod} \ 5).
\label{3-13}
\end{align}
By \eqref{R-3}, we see that
\begin{align}
f_5\equiv f_1^5 \quad ({\rm mod }\ 5). \label{3-14}
\end{align}
Substituting \eqref{3-14}  into  \eqref{3-13} gives
\begin{align}
\sum_{n=0}^\infty \bar{p}(20n+15)q^n \equiv&
3\frac{f_2^{101}}{f_1^{48}f_4^{34}} -\frac{f_2^{29}}{f_{4}^{10}}
-2q\frac{f_2^{77}}{f_1^{40}f_{4}^{18}} +2qf_1^8f_2^5f_4^6
+q\frac{f_2^{29}}{f_1^8f_4^2}
\nonumber\\[6pt]
&\   +3q^2 f_2^5f_4^{14} +3q^2\frac{f_2^{53}}{f_1^{32}f_4^2}
-2q^3\frac{f_2^{29}f_{4}^{14}}{f_1^{24}}
+3q^5\frac{f_4^{46}}{f_1^8f_{2}^{19}} \quad ({\rm mod }\ 5).
\label{3-15}
\end{align}
Combining $2$-dissection formulas \eqref{2-3}, \eqref{2-4} and
congruence \eqref{3-15}, we see that
\begin{align}
\sum_{n=0}^\infty \bar{p}(20n+15)q^n \equiv& 3\frac{f_2^{101}}{
f_4^{34}}\left(
\frac{f_4^{14}}{f_2^{14}f_8^4}+4q\frac{f_4^2f_8^4}{f_2^{10}}
\right)^{12} -\frac{f_2^{29}}{f_{4}^{10}} -2q\frac{f_2^{77}}{
f_{4}^{18}}\left(
\frac{f_4^{14}}{f_2^{14}f_8^4}+4q\frac{f_4^2f_8^4}{f_2^{10}}
\right)^{10}\nonumber\\[6pt]
& \ +2q f_2^5f_4^6
\left(\frac{f_4^{10}}{f_2^2f_8^4}-4q\frac{f_2^2f_8^4}{f_4^2}
\right)^2 +q\frac{f_2^{29}}{ f_4^2}\left(
\frac{f_4^{14}}{f_2^{14}f_8^4}+4q\frac{f_4^2f_8^4}{f_2^{10}}
\right)^{2}
\nonumber\\[6pt]
& \ +3q^2 f_2^5f_4^{14}+3q^2\frac{f_2^{53}}{ f_4^2}\left(
\frac{f_4^{14}}{f_2^{14}f_8^4}+4q\frac{f_4^2f_8^4}{f_2^{10}}
\right)^{8} \nonumber\\[6pt]
&\ -2q^3 f_2^{29}f_{4}^{14} \left(
\frac{f_4^{14}}{f_2^{14}f_8^4}+4q\frac{f_4^2f_8^4}{f_2^{10}}
\right)^{6} +3q^5\frac{f_4^{46}}{f_{2}^{19}}\left(
\frac{f_4^{14}}{f_2^{14}f_8^4}+4q\frac{f_4^2f_8^4}{f_2^{10}}
\right)^{2}\nonumber\\[6pt]
 \equiv& 3\frac{f_4^{134}}{f_2^{67}f_8^{48}}
 -\frac{f_2^{29}}{f_4^{10}}
+2q\frac{f_4^{122}}{f_2^{63}f_8^{40}} +3q\frac{f_2f_4^{26}}{f_8^8}
+q^2\frac{f_4^{110}}{f_2^{59}f_8^{32}}
\nonumber\\[7pt]
& \ +4q^3\frac{f_4^{98}}{f_2^{55}f_8^{24}}
+3q^3f_2^9f_4^2f_8^8+q^4\frac{f_4^{86}}{f_2^{51}f_8^{16}}
+4q^5\frac{f_4^{74}}{f_2^{47}f_8^8}
\nonumber\\[6pt]
& \  +4q^7\frac{f_4^{50}f_8^8}{f_2^{39}}
+q^8\frac{f_4^{38}f_8^{16}}{f_2^{35}}
+4q^9\frac{f_4^{26}f_8^{24}}{f_2^{31}}
+q^{10}\frac{f_4^{14}f_8^{32}}{f_2^{27}}
\nonumber\\[6pt]
&\ + 2q^{11} \frac{f_4^2f_8^{40}}{f_2^{23}}
+3q^{12}\frac{f_8^{48}}{f_2^{19}f_4^{10}} \quad ({\rm mod}\ 5).
\label{3-16}
\end{align}
Extracting the terms with odd powers of $q$ on both sides of
\eqref{3-16}, then dividing by $q$ and replacing $q^2$ by $q$,
 we  reach \eqref{R-1}. This completes the proof. \qed

\section{Proof of Conjecture \ref{conjecture-1}}

In this section, we use
   the $(p, k)$-parametrization of theta
functions given by
 Alaca, Alaca and Williams \cite{Alaca,Alaca-1,Williams}
 to  represent the generating function
  of  $\bar{p}(40+35) $ modulo 5 as a linear combination of several
    functions in
  $p$ and $k$, where $p$ and $k$ are defined in terms of theta
  function  $\varphi(q)$ as given by
  \begin{align}
p =\frac{\varphi^2(q)-\varphi^2(q^3)}{2\varphi^2(q^3)} \label{2-7}
\end{align}
and \begin{align} k =\frac{\varphi^3(q^3)}{\varphi(q)},\label{2-8}
\end{align}  see  Alaca,  Alaca and   Williams \cite{Alaca}.
 Williams \cite{Williams} proved that
\begin{align}
p=2\frac{f_2^3f_3^3f_{12}^6}{f_1f_4^2f_9^6}. \label{R-8}
\end{align}
It turns out that the coefficients of linear combination
  are divisible by $5$. This confirms the conjecture of
  Hirschhorn and Sellers.
More precisely, we have the following congruence.

%

  \begin{theorem}\label{Th-2}  For any nonnegative integer $n$, we have
  \begin{align}
\bar{p}(40n+35)\equiv 0 \quad ({\rm mod}\  5). \label{2-5}
\end{align}

  \end{theorem}

\noindent{\it Proof.} The following representations  of
$q^{\frac{1}{24}}f_1, \ q^{\frac{1}{12}}f_2$ and
 $q^{\frac{1}{6}}f_4$ in terms of $p$ and
$k$ are due to  Alaca and Williams  \cite{Alaca-1},
\begin{align}
q^{\frac{1}{24}}f_1&
=2^{-\frac{1}{6}}p^{\frac{1}{24}}(1-p)^{\frac{1}{2}}(1 +
p)^{\frac{1}{6}}(1 + 2p)^{\frac{1}{8}}
(2 + p)^{\frac{1}{8}}k^{\frac{1}{2}},\label{2-11}\\[6pt]
 q^{\frac{1}{12}}f_2& =2^{-\frac{1}{3}} p^{\frac{1}{12}}(1-p)^{\frac{1}{4}}(1 +
p)^{\frac{1}{12}}(1 + 2p)^{\frac{1}{4}}(2 + p)^{\frac{1}{4}}
k^{\frac{1}{2}}\label{2-12}
\end{align}
and
\begin{align}
q^{\frac{1}{6}}f_4& =2^{-2/3}p^{\frac{1}{6}}(1-p)^{\frac{1}{8}}(1 +
p)^{\frac{1}{24}}(1 + 2p)^{\frac{1}{8}} (2 +
p)^{\frac{1}{2}}k^{\frac{1}{2}}.\label{2-13}
\end{align}
Substituting   \eqref{2-11}, \eqref{2-12} and \eqref{2-13}
 into \eqref{R-1},  we find that
\begin{align}
2^{19}\sum_{n=0}^\infty \bar{p}(40n+35) q^n\equiv
\frac{\sqrt{2}p^{7/8}(1+2p)^{21/8}(2+p)^{21/8}}
{16q^{7/8}(1-p)^6(1+p)^2\sqrt{k}} F(p,k)\quad  ({\rm mod}\
5),\label{R-5}
\end{align}
where $F(p,k)$ is defined by
\begin{align}
F(p,k)=&2621440+30146560p +443678720p^2+4203806720p^3+25364889600p^4
\nonumber\\[6pt]
&  \quad  +112351805440p^5
+378957086720p^6   +980173332480p^7 \nonumber\\[6pt]
& \quad +1961928110080p^8 +3051430471680p^9
 +3658168560640p^{10}
 \nonumber\\[6pt]
&\quad   +3316049272320p^{11}
+2205104730880p^{12}+1020945279360p^{13}
\nonumber\\[6pt]
& \quad   +295430818880p^{14}+40648474720p^{15} +694662000p^{16}
\nonumber\\[6pt]
&  \quad   +12386590p^{19} -82928860p^{18}+168540920 p^{17}+98305
p^{20} .\label{2-14}
\end{align}
By  \eqref{2-11} and \eqref{2-12}, we have
\begin{align}
\frac{f_{2}^{22}}{f_1^{23}}=\frac{\sqrt{2}p^{7/8}(1+2p)^{21/8}(2+p)^{21/8}}
{16q^{7/8}(1-p)^6(1+p)^2\sqrt{k}}. \label{R-6}
\end{align}
Hence  \eqref{R-5} can be rewritten as
\begin{align}
2^{19}\sum_{n=0}^\infty \bar{p}(40n+35) q^n\equiv
 \frac{f_2^{22}}{f_{1}^{23}} F(p,k)\quad  ({\rm mod}\
5),\label{R-7}
\end{align}
where $F(p,k)$ are defined by  \eqref{2-14}. Clearly,
$\frac{f_2^{22}}{f_1^{23}}$ is a formal power series in $q$ with
integer coefficients. By  \eqref{2-8} and \eqref{R-8}, we see that
$p$ and $k$ are also formal power series
 in $q$ with integer coefficients.
It can be seen that the coefficients of $F(p,k)$ are divisible by
$5$. So we reach  the assertion that for all $n\geq 0$,
   $
\bar{p}(40n+35)\equiv 0\ ( {\rm mod }\ 5) $ for $n\geq 0$. This
completes the proof. \qed

To complete the proof of  Conjecture \ref{conjecture-1}, we recall
that Hirschhorn and Sellers \cite{Hirschhorn-1},
 and  Fortin,   Jacob and   Mathieu \cite{Fortin} independently
  derived the congruence
\begin{align}
\bar{p}(4n+3)\equiv 0\quad ({\rm mod}\ 8), \label{3-19}
\end{align}
   for   $n\geq 0$.
This yields
\begin{align}
\bar{p}(40n+35) \equiv 0 \quad ({\rm mod}\ 8),\label{3-20}
\end{align}
for $n\geq 0$. Combining the above congruence \eqref{3-20} and the
congruence $\bar{p}(40n+35) \equiv 0 \ ({\rm mod}\ 5)$ for $n\geq
0$, we come to the conclusion that $ \bar{p}(40n+35) \equiv 0 \
({\rm mod}\ 40)$ for $n\geq 0$.

  \vspace{1cm}
 \noindent{\bf Acknowledgments.}
  This work was supported by the 973 Project, the PCSIRT Project of the
Ministry of Education, and the National Science Foundation of China.

\end{document}